\documentclass{amsart}

\usepackage{amssymb,latexsym}

\newtheorem{theorem}{Theorem}
\newcommand{\bt}{\begin{theorem}}
\newcommand{\et}{\end{theorem}}
\newtheorem{lemma}{Lemma}
\newcommand{\bl}{\begin{lemma}}
\newcommand{\el}{\end{lemma}}
\newtheorem{corollary}{Corollary}
\newcommand{\bc}{\begin{corollary}}
\newcommand{\ec}{\end{corollary}}
\newtheorem{problem}{Problem}
\newcommand{\bprob}{\begin{problem}}
\newcommand{\eprob}{\end{problem}}
\newcommand{\beq}{\begin{equation}}
\newcommand{\eeq}{\end{equation}}
\newcommand{\benum}{\begin{enumerate}}
\newcommand{\eenum}{\end{enumerate}}
\newcommand{\N}{\ensuremath{ \mathbf N }}

\newcommand{\R}{\ensuremath{\mathbf R}}
\newcommand{\Rn}{\ensuremath{\mathbf R^n}}

\newcommand{\mbF}{\ensuremath{ \mathbf F}}

\newcommand{\mcb}{\ensuremath{ \mathcal B}}

\newcommand{\mcf}{\ensuremath{ \mathcal F}}

\newcommand{\mck}{\ensuremath{ \mathcal K}}
\newcommand{\mcl}{\ensuremath{ \mathcal L}}

\newcommand{\mcp}{\ensuremath{ \mathcal P}}

\newcommand{\mcs}{\ensuremath{ \mathcal S}}

\newcommand{\mba}{\ensuremath{ \mathbf a}}

\newcommand{\mbe}{\ensuremath{ \mathbf e}}

\newcommand{\mbo}{\ensuremath{ \mathbf 0}}
\newcommand{\mbu}{\ensuremath{ \mathbf u}}
\newcommand{\mbv}{\ensuremath{ \mathbf v}}

\newcommand{\mbx}{\ensuremath{ \mathbf x}}

\DeclareMathOperator{\conv}{\text{conv}}

\newcommand{\bmat}{\left(\begin{matrix}}
\newcommand{\emat}{\end{matrix}\right)}
\newcommand{\bsmallmat}{\left(\begin{smallmatrix}}
\newcommand{\esmallmat}{\end{smallmatrix}\right)}

\DeclareMathOperator{\qqand}{\qquad\text{and}\qquad}

\title[Finitely many  implies infinitely many]{Finitely many   implies infinitely many}
\author{Melvyn B. Nathanson}
\address{Department of Mathematics\\Lehman College (CUNY)\\
Bronx, NY 10468}
\email{melvyn.nathanson@lehman.cuny.edu}
\date{\today}



\thanks{Supported in part by PSC CUNY Grant \# 66197-00 54.}

\begin{document}

\begin{abstract}
Many mathematical statements have the following form.      
 If something is true for all finite subsets of an infinite set $I$, 
 then it is true for  all of $I$.  
This paper describes some old and new results on infinite sets of linear 
and polynomial equations with the
property that solutions for all finite subsets of the set of equations implies 
the existence of a solution for the infinite set of equations. 
\end{abstract}

\maketitle

\section{Introduction}
Many mathematical statements have the following form.      
 If something is true for all finite subsets of an infinite set $I$, 
 then it is true for  all of $I$.  
 
 Here are two classical examples.\footnote{For completeness, the Appendix 
 contains proofs of the de Bruijn-Erd\H os theorem and  Helly's theorem.}      
 A graph $G$ with vertex set $V$ and edge set $E$ 
 is \emph{$k$-colorable} if there is a function $\chi:V \rightarrow \{1,2,\ldots, k\}$ such that, 
 if $\{v, v'\}$ is an edge in $E$, then $\chi(v) \neq \chi(v')$. Let $G$ be a graph with infinitely many vertices.  The de Bruijn-Erd\H os theorem~\cite{debr-erdo51}  
 states that if every finite subgraph of $G$ is $k$-colorable, 
 then the infinite graph $G$ is $k$-colorable.  
  
Let $\mck = \{K_i:i \in I\}$ be an infinite set of closed convex subsets of \Rn\ that contains at least one compact 
convex set.  Helly's theorem~\cite{hell23} states that if every finite subset of \mck, 
or, simply, every $n+1$ sets in \mck, have a nonempty intersection, 
then the infinite intersection $\bigcap_{i \in I} K_i$ is nonempty.

We shall  describe some old and new results on infinite sets of linear 
and polynomial equations with the
property that solutions for all finite subsets of the set of equations implies 
the existence of a solution for the infinite set of equations. 
These are the notes of a talk in the New York Number Theory Seminar on 
November 30, 2023.

\section{Equations over finite sets}
 Let $J$ be a nonempty set.  Let $X = \{x_j: j \in J\}$ be a set of variables and 
let $\Omega = \prod_{j \in J} \Omega_j$ be the Cartesian product of a family 
of nonempty sets.  We denote by $f(X)$ a statement about the variables 
$X = \{x_j: j \in J\}$.  
Other, equally vague words, such as ``sentence'' or ``proposition'' or ``assertion'', 
might substitute for ``statement''.
Let $\Lambda$ be the set of $J$-tuples $U = (u_j)_{j \in J} \in \Omega$ 
for which the statement $f(U)$ is true.

For example, the statement could be an equation of the form 
\[
f(X) = 0
\] 
in the variables $X = \{x_j: j \in J\}$, and the 
$J$-tuples $U = (u_j)_{j \in J} \in \Omega$ for which the statement is true
are the solutions of the equation, that is, the $J$-tuples $U \in \Omega$ 
such that $f(U) = 0$. 

Here is a simple result in the case of finite sets $\Omega_j$. We recall  the 
\emph{finite intersection property} of subsets of a compact set: If $\Lambda = \{\Lambda_i:i \in I\}$ 
is an infinite set of closed subsets of a compact space and if the intersection 
of every finite subset of $\Lambda$ is nonempty, then the infinite intersection 
$\bigcap_{i\in I} \Lambda_i$ is nonempty. 

\bt
Let $J$ be a nonempty finite or infinite set. 
Let $X = \{x_j:j \in J\}$ be a set of variables 
and let $\Omega = \prod_{j \in J} \Omega_j$ be the Cartesian product of a family 
of nonempty finite sets $\Omega_j$.  
Let $I$ be an infinite set, and, for all $i \in I$, let  $f_i(X)$ 
be a statement in which only finitely many variables $x_j \in X$ occur. 
The infinite set of statements $ \{f_i(X) : i\in I\}$ 
has a solution in $\Omega$ if and only if, for every finite subset $S$ of $I$, 
the set of statements $ \{f_i(X) : i\in S\}$  has a solution in $\Omega$.
\et
   
\begin{proof}
If the infinite set   of statements $ \{f_i(X) : i\in I\}$ has a solution, then every finite subset 
of $ \{f_i(X) : i\in I\}$ has a solution. 

Conversely, suppose that  every finite subset of $S$ of $I$, 
the set of statements $ \{f_i(X) : i\in S\}$  has a solution. 
With the discrete topology, the finite sets $\Omega_j$ are compact 
and so $\Omega$ is compact 
with the usual product topology (Tikhonov's theorem).  
For each $i \in I$, let $\Lambda_i$ be the set of solutions of the 
statement $f_i(X)$ in $\Omega$.  Because the statement $f_i(X)$  
contains only finitely many variables $x_j \in X$ and because the sets $\Omega_j$ 
are finite, the solution set $\Lambda_i$ is closed in $\Omega$.  
Because every finite set of statements has a solution, the intersection 
of a every finite set of the sets $\Lambda_i$ is nonempty.  
The finite intersection property  
implies that $\bigcap_{i\in I} \Lambda_i$ is nonempty, and so the infinite 
set $ \{f_i(X) : i\in I\}$  of statements has a solution in $\Omega$. 
This completes the proof. 
\end{proof}

Let $J$ be a nonempty finite or infinite set and let $\mbx = \{x_j:j\in J\}$ be a set of variables.
We denote by $f(\mbx)= 0$ an equation in finitely or infinitely many of the variables $x_j \in \mbx$.  
Let $I$ be an infinite set and let $ \{f_i(\mbx) = 0: i \in I\}$ be a set of equations. 
Then $\mbu = (u_j)_{j\in J}$ is a \emph{common solution}, or, simply, a \emph{solution} 
of the equations in $ \{f_i(\mbx) = 0: i \in I\}$  if 
$f_i(\mbu) = 0$ for all $i \in I$. 

\bc[Abian~\cite{abia72,abia73}]
Let $R$ be a finite ring, not necessarily commutative. 
An infinite set of polynomial equations in infinitely many 
not necessarily commutative variables 
over $R$ has a  solution in $R$ if and only if every finite subset 
of the set of equations has a solution in $R$. 
\ec

\section{Quadratic equations that do not satisfy ``finitely many implies infinitely many''} 

It is an important observation that there are sets of polynomial equations 
such that every proper subset 
of the set of equations has a common solution 
but not the entire set. Here is a quadratic example due to Alexander Abian~\cite{abia76}.  

Let \mbF\ be a field, let $J = \mbF \cup\{*\}$, and let $\mbx = \{ x_* \} \cup \{x_a:a\in \mbF\}$
be a set of variables.  Consider the set of quadratic equations 
\[
Q_{\mbF} =  \{(a-x_*)x_a-1=0: a \in \mbF\}.
\]
If $ \mbu \in \mbF \times \mbF^{\mbF}$ is a common solution
of these equations with $x_* = u \in \mbF$, then, for the equation $(u-x_*)x_u-1=0$, we obtain 
\[
(u-u)x_u-1=0
\] 
which is absurd.  Therefore, the set of equations $Q_{\mbF}$ has no common solution. 
 
Let $b \in \mbF$ and consider the set of equations obtained by deleting  from $Q_{\mbF}$ 
the single equation $(b-x_*)x_b-1=0$.  
The set of equations 
\[
Q_{\mbF}  \setminus \{ (b-x_*)x_b-1=0 \} = \{(a-x_*)x_a -1=0: a \in \mbF \text{ and } a \neq b\}
\]
has the common solution $x_*=b$ and $x_a = \frac{1}{a-b}$  for all $a \neq b$.
In particular, if \mbF\ is an infinite field, then 
every finite subset of the equations in $Q_{\mbF}$ has a common solution 
but the infinite set $Q_{\mbF}$ 
has no common solution. 
Thus, ``finitely many implies infinitely many'' does not necessarily hold for sets of polynomials.

We shall consider three classes of linear equations where the 
``finitely many implies infinitely many'' principle does hold. 
First, we consider  linear equations in finitely many variables, that is, expressions of the form 
\beq                                       \label{InfSys:eqn-1}
\sum_{j=1}^n a_jx_j + b = 0. 
\eeq
Next,  we consider  linear equations in  infinitely many variables, that is, 
expressions  of the form 
\beq                                       \label{InfSys:eqn-2}
\sum_{j \in J} a_jx_j + b = 0
\eeq
where $J$ is an infinite set and $a_j$ is nonzero for only finitely many $j \in J$.
Finally,  we  discuss linear equations in countably infinitely many variables, that is, 
expressions  of the form
\beq                                       \label{InfSys:eqn-3}
\sum_{j=1}^{\infty} a_jx_j + b = 0. 
\eeq
in which there is no restriction on the number of nonzero coefficients $a_j$.

\section{Infinite sets of linear equations in $n$ variables}

\bt
Let $\mcl$ be the set  of linear equations in $n$ variables $\mbx = \{x_1,\ldots, x_n\}$ 
with coefficients in the field \mbF. 
Let $I$ be an infinite set and let $\mcl (I) = \{ \ell_i: i \in I\}$ be a subset of \mcl. 
If every finite subset of $\mcl (I) $ has a common solution, then the infinite set $\mcl(I)$ 
has a common solution. 
\et

We shall prove the following stronger statement.  

Let $\delta_{i,j}$ be the Kronecker delta.

\bt              \label{InfSys:theorem:n+1}
Let  $\mcl$ be the set of linear equations in $n$ variables 
$\mbx = \{x_1,\ldots, x_n \}$  with coefficients in the field \mbF. 
Let $I$ be a nonempty set and let $ \mcl (I)   = \{ \ell_i: i \in I\}$ be a subset of \mcl. 
If every subset  of $\mcl (I) $ of cardinality at most $n+1$ 
has a common solution, then the $\mcl (I) $ has a common solution. 
\et

\begin{proof}
We define   addition and scalar multiplication of linear equations as follows.  
If $\ell_i$ is the equation $\sum_{j=1}^n a_{i,j}x_j + b_i = 0$ for $i=1,2$ 
and if $c_1, c_2 \in \mbF$, then $c_1\ell_1+c_2\ell_2$ is the linear equation 
\[
\sum_{j=1}^n ( c_1 a_{1,j} +  c_2 a_{2,j}) x_j + (b_1+  b_2) = 0.
\]
With these operations, the set \mcl\ is a vector space over the field \mbF\ 
of dimension $n+1$ with basis 
\[
\{\mbe_0,\mbe_1,\ldots, \mbe_n\}
\]
where $\mbe_i$ is the equation 
\[
\sum_{j=1}^n \delta_{i,j} x_j  + 0 = 0 
\]
for  $i \in \{ 1,\ldots, n \}$  and  $\mbe_0$ is the equation
\[
\sum_{j=1}^n 0 x_j + 1= 0.
\]
The linear equation $\sum_{j=1}^n a_i x_i + b = 0$ is the unique  linear combination 
\[
\sum_{j=1}^n a_i \mbe_j +  b\mbe_0. 
\]

Let $W$ be the subspace of $\mcl$ generated by $\mcl (I) $ 
and let $\dim(W) = m$.  We have $0 \leq m \leq n+1$.  
Let $\mcb = \{\ell_1,\ldots, \ell_m\}$ be a set of linear equations that is a basis for $W$, 
and let $\mbu = (u_1,\ldots, u_n)$ be a common solution of the $m$ equations in \mcb.
Every equation in $S$ is a linear combination of the $m$ equations in \mcb, 
and so  $\mbu = (u_1,\ldots, u_n)$ is a common solution of the  equations in $S$.
This completes the proof. 
\end{proof}

Let $\mbF$ be an algebraically closed field.  
Hilbert's Nullstellensatz says that every set $\Phi$ of polynomials in $\mbF[x_1,\ldots, x_n]$ 
has a common solution in $\mbF$ if and only if the ideal generated by $\Phi$ 
is a proper ideal of $\mbF[x_1,\ldots, x_n]$.     
This is equivalent to the statement that $\Phi$ has a common solution 
if and only if every finite subset of $\Phi$ has a common solution 
(because, if every finite subset of $\Phi$ has a common solution, 
then the equation $\sum_{i=1}^m f_ig_i=1$ with polynomials 
$f_i \in \Phi$ and $g_i \in \mbF[x_1,\ldots, x_n]$ is impossible.  

We can modify the proof of Theorem~\ref{InfSys:theorem:n+1} 
to prove a version of the Nullstellensatz for polynomials 
with coefficients in a Noetherian ring. 

Let $\Phi$ be a nonempty set of polynomials.  
If every finite subset of $\Phi$ has a common solution, then 
the set $\Phi$ generates a proper ideal of $R[x_1,\ldots, x_n]$ 
(because the ideal cannot contain 1) .

\bt
Let $R$ be a Noetherian ring.  If $\Phi$ is a set of polynomials in $R[x_1,\ldots, x_n]$ such 
that every finite subset of $\Phi$ has a common solution, then $\Phi$ has a common solution. 
\et

\begin{proof}
By Hilbert's basis theorem, the polynomial ring $R[x_1,\ldots, x_n]$ is Noetherian.  
Let $\mcl( \Phi )$ be the ideal generated by  $\Phi$, and let $\{\mbe_1, \ldots, \mbe_m\}$ 
be a finite set of polynomials that generates $\mcl( \Phi )$. 
For all $i \in \{ 1,\ldots, m\}$, there is a finite set of polynomials $\{f_{i,1},\ldots, f_{i,k_i}\}$ 
in $\Phi$ and a finite set of polynomials $\{g_{i,1},\ldots, g_{i,k_i}\}$ in $\R[x_1,\ldots, x_n]$ 
such that
\[
\mbe_i = \sum_{r_i=1}^{k_i} f_{i,r_i} g_{i,r_i}. 
\]
The finite set of polynomials 
\[
\mcf = \{f_{i,r_i}:i \in \{ 1,\ldots, m\} \text{ and } r_i \in \{ 1,\ldots, k_i \} \} 
\]
has a common solution $\mbu = (u_1,\ldots, u_n) \in R^n$, and so $\mbu = (u_1,\ldots, u_n)$ 
is a common solution of the polynomials in the generating set $\{\mbe_1, \ldots, \mbe_m\}$. 
Every polynomial in $\Phi$, and, indeed, every polynomial in the ideal $\mcl( \Phi )$, 
is of the form $\sum_{i=1}^m h_i \mbe_i$ for polynomials $h_i \in R[x_1,\ldots, x_n]$, 
and so $\mbu = (u_1,\ldots, u_n)$ 
is a common solution of the polynomials in the $\mcl( \Phi )$. 
This completes the proof. 
\end{proof}

\section{Infinite sets of linear equations in infinitely many variables, 
but with only finitely many summands}

The following beautiful argument is due to Alexander Abian.  

\bt[Abian~\cite{abia76}] 
Let \mbF\ be a field and let $I$ and $J$ be infinite sets.  
Let $\mbx = (x_j)_{j\in J}$ be a sequence of variables 
and, for all $i \in I$,  let 
\[
f_i(\mbx) = \sum_{j \in J} a_{i,j}x_j + b_i = 0
\]
be a linear equation with coefficients $a_{i,j} \in \mbF$, $b_i \in \mbF$, 
and $a_{i,j} \neq 0$ for only finitely many $j$.  
Consider the infinite set of equations $\mcl(I) = \{f_i(\mbx) = 0: i \in I\}$.
If, for every finite subset $I_0$ of $I$, the finite set of equations 
$\mcl(I_0)  =  \{f_i(\mbx) = 0 :i \in I_0\}$ has a common solution, 
then the infinite set $\mcl(I) $  has a common solution. 
\et

\begin{proof}
Let $\mcl$ be the infinite-dimensional vector space 
whose vectors are linear equations of the form 
\[
f(\mbx) =  \sum_{j \in J} a_jx_j + b = 0
\]
where $a_j \in \mbF$, $b \in \mbF$, and $a_j \neq 0$ for only finitely many $j \in J$. 
Let  
\[
\mcl(I) = \left\{ f_i(\mbx)= 0 : i \in I\right\} 
= \left\{ \sum_{j \in J} a_{i,j}x_j + b_i = 0 : i \in I\right\} 
\]
and let $W$ be the subspace of $\mcl$ generated by $\mcl(I)$. 

Let $\mbe_0 \in \mcl$ be the equation 
\[
\sum_{j\in J} 0 x_j + 1 = 0.
\]
We shall prove that $\mbe_0 \notin W$. 

If $\mbe_0 \in W$, then there is a finite subset $I_0 $ of $I$, 
a finite set of equations  
\[
\mcl(I_0) = \left\{ \sum_{j \in J} a_{i,j}x_j + b_{i} = 0 : i \in I_0  \right\} \subseteq \mcl(I) 
\] 
and a finite set of scalars 
$\{\lambda_i : i \in I_0 \}$ such that 
\[
\mbe_0 = \sum_{i \in I_0} \lambda_i \left( \sum_{j \in J} a_{i,j}x_j + b_{i} = 0 \right). 
\]
It follows that the linear equations 
\[
 \sum_{j\in J} 0 x_j  + 1 = 0 
\]
and 
\[ 
  \sum_{i \in I_0} \lambda_i \left( \sum_{j \in J} a_{i,j}x_j + b_{i} = 0 \right)
\]
are the same equation.  
Equivalently, 
\[
 \sum_{j\in J} 0 x_j  + 1 = 0 
\]
and 
\[
\sum_{j \in J} \left(  \sum_{i \in I_0} \lambda_i a_{i,j} \right) x_j  
 +  \sum_{i \in I_0} \lambda_i b_{i}  
= 0 
\]
are the same equation. 
Because the finite set of equations  $\mcl(I_0)$ has a common solution, 
there exist scalars $\mbu = (u_j)_{j \in J}$ such that 
\[
 \sum_{j \in J} a_{i,j}u_j + b_{i}  = 0 \qquad\text{for all $i \in I_0$} 
 \]
and so 
\[
\sum_{j \in J} \left(  \sum_{i \in I_0} \lambda_i a_{i,j} \right) u_j  
 +  \sum_{i \in I_0} \lambda_i b_{i}  
= 0 
\]
It follows that 
\[
0 =  \sum_{j\in J} 0 u_j  = 1
\]
which is absurd.  Therefore, $\mbe_0 \notin W$.

Let $\mcb_W$ be a basis for the subspace $W$ 
generated by $\mcl(I)$  in the vector space $\mcl$.  
Then $\mcb_W \cup \{\mbe_0\}$ is a basis for the subspace 
$W'$ generated by $\mcl(I) \cup \{\mbe_0\}$.  
We can extend the set $\mcb_W \cup \{\mbe_0\}$ to a basis $\mcb$ 
for the vector space $\mcl$. 
Let $\alpha:\mcl \rightarrow \mbF$ be the unique linear functional defined by 
\[
\alpha(\mbv) = \begin{cases}
0 & \text{if $\mbv \in \mcb \setminus \{\mbe_0 \}$} \\
1& \text{if $\mbv = \mbe_0$.}
\end{cases}
\]
It follows that the restriction of $\alpha$ to $W$ is 0, 
and so the restriction of $\alpha$ to the set $\mcl(I)$ is 0. 

For all $j \in J$, let $\mbe_j \in \mcl$ be the linear equation $x_j = 0$ 
and let    
\[
u_j = \alpha(\mbe_j) \in \mbF.
\]
We shall prove that $\mbu = (u_j)_{j\in J}$ is a common solution of the equations in $\mcl(I)$.

For all  $i \in I$, the linear equation $\sum_{j\in J} a_{i,j}x_j + b_i = 0$ in $\mcl(I)$ 
can be written in the form 
\[
\sum_{j\in J} a_{i,j} \mbe_j + b_i \mbe_0 
\]
and so  
\begin{align*}
0 & = \alpha\left( \sum_{j\in J} a_{i,j}x_j + b_i = 0 \right) \\
& = \alpha\left(\sum_{j\in J} a_{i,j} \mbe_j + b_i \mbe_0 \right) \\
& = \sum_{j\in J} a_{i,j}\alpha\left( \mbe_j  \right)+ b_i \alpha\left(\mbe_0 \right) \\
& = \sum_{j\in J} a_{i,j} u_j+ b_i \alpha.
\end{align*}
This completes the proof. 
\end{proof}

\section{Infinite sets of linear equations in infinitely many variables 
and with infinitely many summands}

Let $p$  and $q$ satisfy $p > 1$ and 
\[
\frac{1}{p} + \frac{1}{q} = 1. 
\]
For sequences $\mba = (a_j)_{j=1}^{\infty}  \in \ell^p$ and 
$\mbu = (u_j)_{j=1}^{\infty} \in \ell^q$, we define the sequence 
\[
(\mba,\mbu) = (a_j u_j)_{j=1}^{\infty}. 
\] 
By H\" older's inequality, $(\mba,\mbu)  \in \ell^1$ 
and 
\[
\| (\mba,\mbu) \|_1 = \sum_{j=1}^{\infty} |a_j u_j|    \leq \|\mba\|_p \|\mbu\|_q.
\]
Let $I$ be an infinite set, and let 
 $\mba_i = (a_{i,j})_{j=1}^{\infty} \in \ell^p$ and $b_i \in \R$ for all $i \in I$. 
For all $i \in I$ we have the linear equation in infinitely many variables 
\beq          \label{InfSys:AE}
 (\mba_i,\mbx)  =  \sum_{j=1}^{\infty} a_{i,j} x_j = b_i. 
\eeq
Consider the set $\mcl(I)$ of infinitely many  linear equations in infinitely many variables 
\[
\mcl(I) = \left\{ (\mba_i,\mbx)  = b_i : i \in I\right\} 
= \left\{ \sum_{j=1}^{\infty} a_{i,j} x_j = b_i : i \in I\right\}. 
\]
A solution of the set $\mcl(I)$ of equations is a sequence $\mbu \in \ell^q$ such that 
$(\mba_i, \mbu )= b_i$ for all $i \in I$. 

It is \emph{not true} that every  set of linear equations in infinitely many 
variables satisfies the ``finitely many implies infinitely many`` paradigm.  
Here is a simple example of an infinite set of linear equations for which every finite subset 
has a solution but the infinite set has no solution. 
Let $\N = \{1,2,3,\ldots\}$ be the set of positive integers.  
For all $i \in \N$, let $\mba_i = (a_{i,j})_{j=1}^{\infty}$ be the sequence 
defined by 
\[
a_{i,j} = \begin{cases}
1 & \text{if $j \leq i$} \\ 
0 & \text{if $j > i$.} 
\end{cases}
\]
For all $p \geq 1$ we have $\mba_i \in \ell^p$  and 
\[
\| \mba_i \|_p = i^{1/p}.
\] 
Let $b_i = i$ for all $i \in \N$. 
Consider the infinite set of linear equations 
\[
\mcl(\N) = \{(\mba_i,\mbx) = i:i \in \N \}  = \left\{ \sum_{j=1}^i x_j = i:i \in \N \right\}.
\]
For all $n \in \N$, define the sequence $\mbu^{(n)} = (u_j)_{j=1}^{\infty}\in \ell^q$ by 
\[
u_j = \begin{cases}
1 & \text{if $j \leq n$} \\ 
0 & \text{if $j > n$.} 
\end{cases}
\]
For all $i \in \{1,2,\ldots, n\}$ we have  
\beq                    \label{bbb}
(\mba_i,\mbu^{(n)}) = \sum_{j=1}^i u_j =  i
\eeq
and so every finite subset of the infinite set $\mcl(\N) $ has a solution in $\ell^q$. 

If a sequence $\mbu \in \ell^q$ satisfies every equation 
in the set $\mcl(\N) $, then  we have 
\[
i =  (\mba_i, \mbu)  = \left\| (\mba_i, \mbu ) \right\|_1 \leq \left\| \mba_i \right\|_p \left\| \mbu\right\|_q 
= i^{1/p} \left\| \mbu \right\|_q 
\]
and so 
\[
i^{1/q} = i^{1-(1/p)} \leq \left\| \mbu \right\|_q < \infty
\]
for all $i \in \N$, which is absurd. 
Thus, the set $\mcl(\N) $ has no solution in $\ell^q$.

With an additional condition, one can prove that 
sets of linear equations in infinitely many variables 
do satisfy ``finitely many implies infinitely many``.  
The following theorem of Abian and Eslami~\cite{abia-esla82}  
extends a classical result of F.   Riesz~\cite{ries13}.

Let $M > 0.$ 
 Extending a classical result of F. Riesz~\cite{ries13} 
 (see also Banach~\cite[p. 47]{bana87}),  
 Abian and Eslami~\cite{abia-esla82} proved that if 
$\mcl(I) = \left\{ (\mba_i, \mbx) = b_i : i \in I \right\}$ is an infinite set 
of linear equations such that every finite subset $\mcl(I_0) $ has a solution 
$\mbu^{(0)} \in \ell^q$ 
with $\| \mbu^{(0)}\|_q \leq M$, then the infinite set of equations has a solution 
$\mbu \in \ell^q$ with $\| \mbu\|_q \leq M$.

\bprob
Is there a condition weaker than that every finite subset of an infinite set $S$ of 
linear equations has a solution 
\mbx\ with $\ell_q$ norm at most $M$ that still guaranties a solution of every equation in $S$?
\eprob

The linear equation~\eqref{InfSys:AE} can be written in the form 
\[
\sum_{j=1}^{\infty} f_{i,j}(x_j) = b_i 
\]
where 
\[
 f_{i,j}(x_j) = a_{i,j}x_j 
\]
is a  linear polynomial  with constant term 0.  
Suppose we replace  the linear polynomial with a polynomial of degree at most $d$
and  with constant term 0:
\[
f_{i,j}(x_j) = \sum_{k=1}^d a_{i,j,k} x_j^k.
\]
For all $i \in \N$ and $b_i \in \R$ we obtain a polynomial equation in infinitely many variables 
\[
\sum_{j=1}^{\infty} f_{i,j}(x_j) = b_i 
\]
Nathanson and Ross~\cite{nath-ross24x} proved the following.

\bt
Let $d$ be a positive integer and let $q > d$.  
For all $i \in \N$ and $k \in \{1,2,\ldots, d\}$, let 
\[
\left(a_{i,j,k}\right)_{j=1}^{\infty} \in \ell^{q/(q-k)}.   
\]
Define the polynomials  
\[
f_{i,j}(x_j) = \sum_{k=1}^d a_{i,j,k} x_j^k 
\]
of degree at most $d$ with constant term 0.
For every sequence of real numbers $(b_i)_{i=1}^{\infty}$, 
consider the countably infinite set \mcs\ of polynomial equations 
 in infinitely many variables 
\[
\mcl(\N) = \left\{ \sum_{j=1}^{\infty} f_{i,j}(x_j) = b_i : i \in \N \right\}
\]
Let $M > 0$.  
There exists a sequence $\mbu  = (u_i)_{i=1}^{\infty}  \in \ell^q$ such that $\|\mbu\|_q \leq M$ 
and $\mbu$ is a common solution  of the equations in $\mcl(\N)$ if and only if 
every finite subset of $\mcl(\N)$ has a common solution $\mbu^{(0)} \in \ell^q$ 
with $\|\mbu^{(0)}\|_q \leq M$.
\et


\appendix
\section{Theorems of Rado and de Bruijn-Erd\H os}

Let $X$ be a set. 
Let  $\mcp(X)$ be the  \emph{power set}\index{power set} of   $X$, 
that is,  the set of all subsets of $X$ and let  $\mcf^*(X)$ 
be the set of all nonempty finite subsets of $X$.  

Let $I$ be a set and let $\{X_i:i \in I\}$ be a sequence of nonempty subsets of $X$.
A \emph{choice function}\index{choice function} is a function 
$\varphi:I \rightarrow \mcp(X)$ such that 
$\varphi(i) \in X_i$ for all $i \in I$. 


The following is \index{Rado's selection lemma}\emph{Rado's selection lemma}.

\bt[Rado~\cite{rado49}]                                 \label{InfSys:theorem:Rado49}
Let $I$ be an infinite set, and, for all $i \in I$,
 let $X_i$ be a nonempty finite set.  
For all $A \in \mcf^*(I)$, let $\varphi_A: A \rightarrow \bigcup_{i\in I} X_i$ be a choice function, 
that is, $\varphi_A(a) \in X_a$ for all $a \in A$. 
There exists a choice function $\Phi:I \rightarrow X$ such that, for all $A \in \mcf^*(I)$, 
there exists $B \in \mcf^*(I)$ 
with $A \subseteq B$ such that  $\varphi_{B}(a) = \Phi(a)$ for all $a \in A$. 
\et

\begin{proof}
The following proof  is due to Gottschalk~\cite{gott51}. 

With the discrete topology, the finite set $X_i$ is compact for all $i \in I$ 
and so the set $\Omega = \prod_{i\in I} X_i$ 
with the usual product topology is compact (Tikhonov's theorem).  
The elements of $\Omega$ are choice functions $\varphi:i \mapsto \Omega$ 
with $\varphi(i) \in X_i$ for all $i \in I$.  

Let $A \in \mcf^*(I)$.   For every set $B \in \mcf^*(I)$ such that $A \subseteq B$, 
the restriction of the choice function $\varphi_B$ to the set $A$ defines a choice function on $A$, that is, the $|A|$-tuple $(\varphi_B(a):a \in A)$ in $\prod_{a\in A} X_a$.

Let $E_A$ be the set of all $\varphi \in \Omega$ such that there exists 
$B \in \mcf^*(I)$ such that $A \subseteq B$ and $\varphi(a) = \varphi_B(a)$ for all $a \in A$.
There exists $\varphi \in \Omega$ such that $\varphi(a) = \varphi_A(a)$ for all $a \in A$, 
and so $E_a$ is nonempty.  
There are only finitely many $|A|$-tuples in $\prod_{a \in A}X_a$ 
that are restrictions of choice functions 
$\varphi_B$ with $A \subseteq B$, and so $E_A$ is a nonempty closed subset of $\Omega$.  
Moreover, if $A,A' \in \mcf^*(I)$ and $A \subseteq A'$, then $\varphi \in E_{A'}$ implies there 
exists $B' \in \mcf^*(I)$ such that $A' \subseteq B'$ and $\varphi(a') = \varphi_{B'}(a')$ 
for all $a' \in A'$.  Because $A\subseteq A'$, we have $\varphi(a) = \varphi_{B'}(a)$ 
for all $a \in A$, and so $\varphi \in E_A$.  Thus, $E_{A'} \subseteq E_A$. 

Let $A_1,\ldots, A_n \in \mcf^*(I)$ and let  $A = \bigcup_{i=1}^n A_i  \in \mcf^*(I)$.  
For all $i \in \{1,\ldots, n\}$, the inclusion $A_i \subseteq A$ implies $E_A \subseteq E_{A_i}$.  
It follows that  $E_A \subseteq \bigcap_{i=1}^n E_{A_i}$ and so the closed set 
$\bigcap_{i=1}^n E_{A_i}$ is nonempty. 
Thus, the set of closed sets $\{E_A:A\in \mcf^*(I)$ has the finite intersection property.  
By compactness, there exists $\Phi \in \bigcap_{A \in \mcf^*(I)} E_A$.
This completes the proof. 
\end{proof}

\bc[Gottschalk~\cite{gott51}] 
Let $\{X_i:i \in I\}$ be an infinite set of nonempty finite sets.  
The set $X = \bigcup_{i\in I} X_i$ has a one-to-one choice function if and only if 
every finite subset has a one-to-one choice function. 
\ec

A \emph{$k$-coloring} of the vertices of a graph $G = (V,E)$ is a function 
$\chi:G \rightarrow \{1,2,\ldots, k\}$ 
such that if $\{v_1,v_2\}$ is an edge, then $\chi(v_1) \neq \chi(v_2)$. 
The graph $G = (V,E)$ is \emph{infinite} if its vertex set $V$ is infinite.

\bt[de Bruijn and Erd\H os~\cite{debr-erdo51}]                \label{InfSys:theorem:deBruijn}
Let $G$ be an infinite graph.  If every finite subgraph of $G$ is $k$-colorable, then $G$ is $k$-colorable.
\et

\begin{proof}
To every vertex $v \in V$ we assign the set $I = \{1,2,\ldots, k\}$.   
For every finite subset $W$ of $V$, there is a $k$-coloring function $\chi_W$,  
and  $\chi_W$ is a choice function. 
Let $\chi:V \rightarrow I$ be the  choice function  obtained from  Rado's selection  lemma.  
Let $e = \{v_1,v_2\}$ 
be an edge in the graph $G$.  Then $e$ is a 2-element subset of $V$ 
and there is a finite set $W$ of $V$ such that $e \subseteq W$ and $\chi(v_1) = \chi_W(v_1)$
and $\chi(v_2) = \chi_W(v_2)$.  Because $\chi_W$ is a $k$-coloring of a finite subgraph, we have 
$\chi(v_1)\neq \chi(v_2$ and so $\chi$ is a $k$-coloring of $G$. 
This completes the proof. 
\end{proof}

\section{Helly's theorem}

\bt[Radon~\cite{rado21}]                        \label{InfSys:theorem:Radon}
Let $A$ be a set of $n+2$ distinct points in $\R^n$.  
There is a partition of $A$ into nonempty pairwise disjoint sets 
\[
A = A_1 \cup A_2 \qqand  A_1 \cap A_2 = \emptyset 
\]
 such that 
 \[
 \conv(A_1) \cap \conv(A_2) \neq \emptyset. 
 \]
\et

\begin{proof}
Let $A = \{\mba_1,\ldots, \mba_{n+2}\}$, where  
\[
\mba_j  = \bmat a_{1,j} \\ a_{2,j} \\ \vdots \\ a_{n,j} \emat \in \Rn
\]
for  all $j \in \{1,2,\ldots, n+2\}$.  
Every set of $n+1$ homogeneous linear equations 
in $n+2$ variables has a nonzero solution. Thus, 
\begin{align}
\sum_{j=1}^{n+2} a_{i,j} x_j & = 0 \qquad \text{for $i \in \{1,2,\ldots, n\}$}         \label{InfSys:Radon-1} \\ 
\sum_{j =1}^{n+2}  x_j & = 0.    \label{InfSys:Radon-2}
\end{align}
 has a nonzero solution $\mbu = \bsmallmat u_1 \\ \vdots \\ u_{n+2}\esmallmat$. 
Let 
\[
A_1 = \{ \mba_j \in A: u_j > 0\} \qqand A_2 = \{\mba_j \in A : u_j \leq 0\}. 
\] 
 Equation~\eqref{InfSys:Radon-2} implies that the vector $\mbu$ 
 has both positive and negative coordinates, and so 
 the sets $A_1$ and $A_2$ are nonempty and partition $A$.  
Moreover, 
\[
s = \sum_{\mba_j \in A_1} u_j  = \sum_{\mba_j \in A_2} (-u_j) > 0.
\] 
Therefore,
\[
\sum_{\mba_1\in A_1} \frac{u_j}{s} = \sum_{\mba_j \in A_2} \frac{-u_j}{s} = 1.
\]
The set of $n$ linear equations~\eqref{InfSys:Radon-1}
is equivalent to the vector equation 
\[
\sum_{j=1}^{n+2} x_j \mba_j = \mbo.
\]
We obtain the vector relation 
\[
\sum_{ \mba_j\in A_1} u_j\mba_j = \sum_{\mba_j\in A_2} (-u_j) \mba_j. 
\]
Dividing by $s$ gives the convex combinations 
\[
\sum_{ \mba_j\in A_1} \frac{u_j}{s} \mba_j 
= \sum_{\mba_j\in A_2} \frac{-u_j}{s} \mba_j \in \conv(A_1) \cap \conv(A_2) 
\]
and so $\conv(A_1) \cap \conv(A_2) \neq \emptyset$. 
This completes the proof. 
\end{proof}

\bt[Helly~\cite{hell23}]                                       \label{InfSys:theorem:Helly}
Let  $m \geq n+2$ and let $K_1,\ldots, K_m$ be convex subsets of $\R^n$.  
If every $n+1$ sets $K_i$ have a nonempty intersection,
 then $\bigcap_{i=1}^m K_i \neq \emptyset$. 
 
Let $I$ be an infinite set and let $\mck = \{K_i:i\in I\}$ be a set 
of closed convex subsets of \Rn\ that contains at least one compact 
convex set.  If every $n+1$ sets in \mck\ have a nonempty intersection, 
then $\bigcap_{i \in I} K_i$ is nonempty. 
\et

\begin{proof}
The proof is by induction on $m$.  Let $m = n+2$. 
By the induction hypothesis, for all $j \in \{1,\ldots, n+2\}$ there is a 
vector 
\[
\mba_j \in \bigcap_{\substack{i=1\\i \neq j}}^{n+2} K_i.
\]
If $\mba_j = \mba_{j'}$ for some $j \neq j'$, then 
\[
\mba_j = \mba_{j'} \in \bigcap_{\substack{i=1\\i \neq j'}}^{n+2} K_i \subseteq K_j 
\] 
and so $\mba_j \in \bigcap_{i=1}^{n+2} K_i \neq \emptyset$ 
and we are done.  

If $\mba_j \neq \mba_{j'}$ for all $j \neq j'$, then $A = \{\mba_1,\ldots, \mba_{n+2}\}$ 
is a set of $n+2$ distinct vectors in \Rn. 
By Radon's theorem (Theorem~\ref{InfSys:theorem:Radon}), 
there is a partition $A = A_1 \cup A_2$ and a vector $\mbv$ such that 
$\mbv \in \conv(A_1) \cap \conv(A_2)$. 
We shall prove that $\mbv \in \bigcap_{i=1}^{n+2} K_i$.

Let $j \in \{1,\ldots, n+2\}$.  
If $\mba_j \in A_1$, then $\mba_i \in K_j$ for all $i \neq j$.  
It follows that $A_2 \subseteq K_j$ and, because the set $K_j$ is convex, we have 
\[
\mbv \in \conv(A_1)\cap\conv(A_2) \subseteq \conv(A_2) \subseteq K_j. 
\]
Similarly, if $\mba_j \in A_2$, then $\mbv \in K_j$.  
Thus, $\mbv \in K_j$ for all $j \in \{1,\ldots, n+2\}$. This completes the proof for $m = n+2$. 

Let $m \geq n+2$ and assume that Helly's theorem holds for $m$ sets.  
Given $m+1$ convex sets $K_1,\ldots, K_{m-1}, K_m, K_{m+1}$, 
we simply apply Helly's theorem to the $m$ convex sets 
$K_1,\ldots, K_{m-1}, K'_m$ with $K'_m = K_m \cap K_{m+1}$, and we are done. 

Let $I$ be an infinite set, and 
let $\mck = \{K_i:i\in I\}$ be a  set of closed convex sets 
such that every $n+1$ sets in \mck\ have a nonempty intersection.  
Let $K_0 \in K$ be compact.  
Then $\mck_0 = \{K_0 \cap K_i: i\in I \setminus \{0\} \}$ is a set 
of nonempty closed subsets of the compact set $K_0$.  
The finite form of Helly's theorem implies that every finite intersection 
of sets in $\mck_0$ is nonempty.  
The finite intersection property for closed compact sets implies that 
\[
\emptyset \neq \bigcap_{ i\in I \setminus \{0\} }(K_0\cap K_i) \subseteq \bigcap_{i \in I} K_i.
\] 
This completes the proof. 
\end{proof}

\def\cprime{$'$} \def\cprime{$'$}
\providecommand{\bysame}{\leavevmode\hbox to3em{\hrulefill}\thinspace}
\providecommand{\MR}{\relax\ifhmode\unskip\space\fi MR }
\providecommand{\MRhref}[2]{%
  \href{http://www.ams.org/mathscinet-getitem?mr=#1}{#2}
}
\providecommand{\href}[2]{#2}

\end{document}